\documentclass{article}
\usepackage{amssymb}
\usepackage{amsfonts}
\usepackage{amsmath}

\def\QDATOPD#1#2#3#4{\displaystyle\left#1\begin{array}{c}
                     #3 \\
                              #4 \\
                            \end{array}\right #2}

\newtheorem{theorem}{Theorem}
\newtheorem{acknowledgement}[theorem]{Acknowledgement}

\newtheorem{corollary}[theorem]{Corollary}

\newenvironment{proof}[1][Proof]{\noindent\textbf{#1.} }{\ \rule{0.5em}{0.5em}}
\include {mak}
\input{tcilatex}

\begin{document}

\title{{\large Connection between ordinary multinomials, generalized
Fibonacci numbers, partial Bell partition polynomials and convolution powers
of discrete uniform distribution}}
\author{Hac\`{e}ne Belbachir \and Sadek Bouroubi \and Abdelkader Khelladi
\and \ \ \ \ \ \ \ \ \ \ \ \ \  \\
USTHB/ Facult\'{e} de Math\'{e}matiques\\
BP 32, El Alia, 16111 Bab Ezzouar, Alger, Algeria.\\
hbelbachir@usthb.dz and hacenebelbachir@gmail.com\\
sbouroubi@usthb.dz and bouroubis@yahoo.fr\\
akhelladi@usthb.dz and khelladi@wissal.dz}
\maketitle

\begin{abstract}
Using an explicit computable expression of ordinary multinomials, we
establish three remarkable connections, with the $q$-generalized Fibonacci
sequence, the exponential partial Bell partition polynomials and the density
of convolution powers of the discrete uniform distribution. Identities and
various combinatorial relations are derived.
\end{abstract}


\section{Introduction}

\noindent Ordinary multinomials are a natural extension of binomial
coefficients, for an appropriate introduction of these numbers see Smith and
Hogatt (1979), Bollinger (1986) and Andrews and Baxter (1987). These
coefficients are defined as follows:

\noindent Let $q\geq 1$ and $L\geq 0$ be integers. For an integer $%
a=0,1,\ldots ,qL,$ the ordinary multinomial $\binom{L}{a}_{q}$ is the
coefficient of the $a$-th term of the following multinomial expansion%
\begin{equation}
\left( 1+x+x^{2}+\cdots +x^{q}\right) ^{L}=\sum\limits_{a\geq 0}\binom{L}{a}%
_{q}x^{a},  \label{1}
\end{equation}%
with $\binom{L}{a}_{1}=\binom{L}{a}$ (being the usual binomial coefficient)
and $\binom{L}{a}_{q}=0$ for $a>qL.$ Using the classical binomial
coefficient, one has%
\begin{equation}
\binom{L}{a}_{q}=\sum_{j_{1}+j_{2}+\cdots +j_{q}=a}\binom{L}{j_{1}}\binom{%
j_{1}}{j_{2}}\cdots \binom{j_{q-1}}{j_{q}}.  \label{2}
\end{equation}

\noindent Some readily well known established properties are

\begin{description}
\item the symmetry relation%
\begin{equation}
\binom{L}{a}_{q}=\binom{L}{qL-a}_{q},  \label{3}
\end{equation}

\item the longitudinal recurrence relation%
\begin{equation}
\binom{L}{a}_{q}=\sum_{m=0}^{q}\binom{L-1}{a-m}_{q},  \label{4}
\end{equation}

\item and the diagonal recurrence relation%
\begin{equation}
\binom{L}{a}_{q}=\sum_{m=0}^{L}\binom{L}{m}\binom{m}{a-m}_{q-1}.  \label{5}
\end{equation}
\end{description}

\noindent These coefficients, as for usual binomial coefficients, are built
trough the Pascal triangle, known as "Pascal Pyramid", see tables: 1, 2 and
3. One can find the first values of the pyramid in {\small SLOANE }\cite{slo}
as{\small \ A027907} for $q=2$, {\small A008287} for $q=3$ and {\small %
A035343} for $q=4$.

\ \ \ \ \ \ \ \ \ \ \ \ \ \ \ \ \ \ \ \ \ \ \

\noindent As an illustration of recurrence relation, we give the triangles
of trinomial, quadrinomial and pentanomial coefficients:

\begin{center}
Table 1: Triangle of \textit{trinomial} coefficients: $\binom{L}{k}_{2}$
\end{center}

\noindent $%
\begin{array}{cccccccccccc}
L\backslash k & \text{\textit{0}} & \text{\textit{1}} & \text{\textit{2}} &
\text{\textit{3}} & \text{\textit{4}} & \text{\textit{5}} & \text{\textit{6}}
& \text{\textit{7}} & \text{\textit{8}} & \text{9} & \text{\textit{10}} \\
\text{\textit{0}} & \text{1} &  &  &  &  &  &  &  &  &  &  \\
\text{\textit{1}} & \text{1} & \text{1} & \text{1} &  &  &  &  &  &  &  &
\\
\text{\textit{2}} & \text{1} & \text{2} & \text{3} & \text{2} & \text{1} &
&  &  &  &  &  \\
\text{\textit{3}} & \text{1} & \text{3} & \text{6} & \text{7} & \text{6} &
\text{3} & \text{1} &  &  &  &  \\
\text{\textit{4}} & \text{1} & \text{4} & \text{10} & \text{16} & \text{19}
& \text{16} & \text{10} & \text{4} & \text{1} &  &  \\
\text{\textit{5}} & \text{1} & \text{5} & \text{15} & \text{30} & \text{45}
& \text{51} & \text{45} & \text{30} & \text{15} & \text{5} & \text{1}%
\end{array}%
$

\ \ \ \ \ \ \ \ \ \ \ \ \ \ \ \ \ \ \ \ \ \ \ \ \ \ \ \ \ \ \ \ \ \ \

\begin{center}
Table 2: Triangle of \textit{quadrinomial} coefficients: $\binom{L}{k}_{3}$
\end{center}

\noindent $%
\begin{array}{cccccccccccccc}
L\backslash k & \text{\textit{0}} & \text{\textit{1}} & \text{\textit{2}} &
\text{\textit{3}} & \text{\textit{4}} & \text{\textit{5}} & \text{\textit{6}}
& \text{7} & \text{\textit{8}} & \text{\textit{9}} & \text{\textit{10}} &
\text{\textit{11}} &  \\
\text{\textit{0}} & \text{1} &  &  &  &  &  &  &  &  &  &  &  &  \\
\text{\textit{1}} & \text{1} & \text{1} & \text{1} & \text{1} &  &  &  &  &
&  &  &  &  \\
\text{\textit{2}} & \text{1} & \text{2} & \text{3} & \text{4} & \text{3} &
\text{2} & \text{1} &  &  &  &  &  &  \\
\text{\textit{3}} & \text{1} & \text{3} & \text{6} & \text{10} & \text{12} &
\text{12} & \text{10} & \text{6} & \text{3} & \text{1} &  &  &  \\
\text{\textit{4}} & \text{1} & \text{4} & \text{10} & \text{20} & \text{31}
& \text{40} & \text{44} & \text{40} & \text{31} & \text{20} & \text{10} &
\text{4} & \text{1}%
\end{array}%
$

\ \ \ \ \ \ \ \ \ \ \ \ \ \ \ \ \ \ \ \ \ \ \ \ \ \ \ \ \ \ \ \ \

\begin{center}
Table 3: Triangle of \textit{pentanomial} coefficients: $\binom{L}{k}_{4}$
\end{center}

\noindent $%
\begin{array}{cccccccccccccccc}
L\backslash k & \text{\textit{0}} & \text{\textit{1}} & \text{\textit{2}} &
\text{\textit{3}} & \text{\textit{4}} & \text{\textit{5}} & \text{\textit{6}}
& \text{\textit{7}} & \text{\textit{8}} & \text{9} & \text{\textit{10}} &
\text{\textit{11}} & \text{\textit{12}} & \text{\textit{13}} &  \\
\text{\textit{0}} & \text{1} &  &  &  &  &  &  &  &  &  &  &  &  &  &  \\
\text{\textit{1}} & \text{1} & \text{1} & \text{1} & \text{1} & \text{1} &
&  &  &  &  &  &  &  &  &  \\
\text{\textit{2}} & \text{1} & \text{2} & \text{3} & \text{4} & \text{5} &
\text{4} & \text{3} & \text{2} & \text{1} &  &  &  &  &  &  \\
\text{\textit{3}} & \text{1} & \text{3} & \text{6} & \text{10} & \text{15} &
\text{18} & \text{19} & \text{{\small 18+}} & \text{{\small 15+}} & \text{%
{\small 10+}} & \text{{\small 6+}} & \text{{\small 3+}} & \text{1} &  &  \\
\text{\textit{4}} & \text{1} & \text{4} & \text{10} & \text{20} & \text{35}
& \text{52} & \text{68} & \text{80} & \text{85} & \text{80} & \text{68} &
\text{{\small =52}} & \text{35} & \text{20} & \cdots%
\end{array}%
$

\ \ \ \ \ \ \ \ \ \ \ \ \ \ \ \ \ \

\noindent Several extensions and commentaries about these numbers have been
investigated in the literature:

\noindent Brondarenko \cite[1993]{bro} gives a combinatorial interpretation
of ordinary multinomials $\binom{L}{a}_{q}$ as \textit{the number of
different ways of distributing "}$a$\textit{" balls among "}$L$\textit{"
cells where each cell contains at most "}$q$\textit{" balls.}

\noindent Using this combinatorial argument, one can easily establish the
following relation%
\begin{eqnarray}
\binom{L}{a}_{q} &=&\sum_{L_{1}+2L_{2}+\cdots +qL_{q}=a}\binom{L}{L_{1}}%
\binom{L-L_{1}}{L_{2}}\cdots \binom{L-L_{1}-\cdots -L_{q-1}}{L_{q}}  \notag
\\
&=&\sum_{L_{1}+2L_{2}+\cdots +qL_{q}=a}\binom{L}{L_{0},L_{1},L_{2}\cdots
,L_{q-1}}.  \label{rrr}
\end{eqnarray}

\noindent For a computational view of the relation (\ref{rrr}) see Bollinger
\cite[1986]{bol}.
\noindent Andrews and Baxter \cite[1987]{and} have considered the
q-analog generalization\ of ordinary multinomials (see also
\cite[1997]{war1} for an exhaustive bibliography). They have defined
the q-multinomial coefficients
as follows%
\begin{equation*}
\QDATOPD[ ] {n}{a}_{q}^{\left( p\right) }=\sum_{j_{1}+j_{2}+\cdots
+j_{q}=a}q^{\sum_{l=1}^{q-1}\left( L-j_{l}\right)
j_{l+1}-\sum_{l=q-p}^{q-1}j_{l+1}}\QDATOPD[ ] {L}{j_{1}}\QDATOPD[ ] {j_{1}}{%
j_{2}}\ldots \QDATOPD[ ] {j_{q-1}}{j_{q}},
\end{equation*}%
where
\begin{equation*}
\QDATOPD[ ] {L}{a}=\QDATOPD[ ] {L}{a}_{q}=\left\{
\begin{array}{ll}
\left( q\right) _{L}/\left( q\right) _{a}\left( q\right) _{L-a} & \text{if}\
0\leq a\leq L \\
0 & \text{otherwise}%
\end{array}%
\right.
\end{equation*}%
is the usual q-binomial coefficient, and where $\left( q\right)
_{k}=\prod_{m=1}^{\infty }\left( 1-q^{m}\right) /\left( 1-q^{k+m}\right) ,$
is called $q$-series. This definition is motivated by the relation (\ref{2}).

\noindent Another extension, the supernomials, has also been considered by
Schilling and Warnaar \cite[1998]{sch}. These coefficients are defined to be
the coefficients of $x^{a}$ in the expression of $\prod_{j=1}^{N}\left(
1+x+\cdots +x^{j}\right) ^{L_{j}}$

\noindent A refinement of the q-multinomial coefficient is also considered
for the trinomial case by Warnaar \cite[2001]{war2}.

\noindent Barry \cite[2006]{bar} gives a generalized Pascal triangle as $%
\binom{n}{k}_{a\left( n\right) }:=\prod_{j=1}^{k}a\left( n-j+1\right)
/a\left( j\right) ,$ where $a\left( n\right) $ is a suitably chosen sequence
of integers.

\noindent Kallas \cite[2006]{kal} and Noe \cite[2006]{noe} give a
generalization of Pascal's triangle by considering the coefficient of $x^{a}$
in the expression of $\left( a_{0}+a_{1}x+\cdots +a_{q}x^{q}\right) ^{L}.$

\noindent The main goal of this paper is to give some connections of the
ordinary multinomials with the generalized Fibonacci sequence, the
exponential Bell polynomials, and the density of convolution powers of
discrete uniform distribution. We will give also some intersting
combinatorial identities

\section{A simple expression of ordinary multinomials}

\noindent If we denote $x_{i}$ the number of balls in a cell, the previous
combinatorial interpretation given by Brondarenko is equivalent to evaluate
the number of solutions of the system%
\begin{equation}
\left\{
\begin{array}{c}
x_{1}+\cdots +x_{L}=a, \\
0\leq x_{1},\cdots ,x_{L}\leq q.%
\end{array}%
\right.  \label{5b}
\end{equation}%
Now, let us consider the system (\ref{5b}). For $t\in \left] -1,1\right[ $,
we have (see also Comtet \cite[Vol.1, p. 92 (pb 16).]{com})%
\begin{equation*}
\sum\limits_{a\geq 0}\binom{L}{a}_{q}t^{a}=\left( 1+t+\cdots +t^{q}\right)
^{L}=\sum\limits_{0\leq x_{1},\cdots ,x_{L}\leq q}t^{x_{1}+\cdots +x_{L}},
\end{equation*}%
and%
\begin{eqnarray*}
\left( 1+t+\cdots +t^{q}\right) ^{L} &=&\left( 1-t^{q+1}\right) ^{L}\left(
1-t\right) ^{-L} \\
&=&\left( \sum_{j=0}^{L}\left( -1\right) ^{j}\binom{L}{j}_{q}t^{j\left(
q+1\right) }\right) \left( \sum_{j\geq 0}\binom{j+L-1}{L-1}t^{j}\right) .
\end{eqnarray*}%
By identification, we obtain the following theorem.

\begin{theorem}
The following identity holds%
\begin{equation}
\binom{L}{a}_{q}=\sum_{j=0}^{\left\lfloor a/\left( q+1\right) \right\rfloor
}\left( -1\right) ^{j}\binom{L}{j}\binom{a-j\left( q+1\right) +L-1}{L-1}.
\label{6b}
\end{equation}
\end{theorem}

\noindent This explicite relation seems to be important since in contrast to
relations (\ref{2}), (\ref{3}) and (\ref{5}), it allows to compute the
ordinary multinomials with one summation symbol.

\noindent In 1711, de Moivre (see \cite[1731]{moi1} or \cite[1756 3rd ed.
p.39]{moi}) solves the system (\ref{5b}) as the right hand side of (\ref{6b}%
).

\begin{corollary}
We have the following identity%
\begin{equation*}
\sum_{j=0}^{\left\lfloor n/2\right\rfloor }\binom{n}{j}\binom{n-j}{j}%
=\sum_{j=0}^{\left\lfloor n/3\right\rfloor }\left( -1\right) ^{j}\binom{n}{j}%
\binom{2n-3j-1}{n-1}.
\end{equation*}
\end{corollary}

\begin{proof}
It suffices to use relation (\ref{4}) in Theorem 1 for $q=2$ and $a=L=n.$
\end{proof}

\section{Generalized Fibonacci sequences}

\noindent Now, let us consider for $q\geq 1$, the "multibonacci" sequence $%
(\Phi _{n}^{\left( q\right) })_{n\geq -q}$ defined by%
\begin{equation*}
\left\{
\begin{array}{l}
\Phi _{-q}^{\left( q\right) }=\cdots =\Phi _{-2}^{\left( q\right) }=\Phi
_{-1}^{\left( q\right) }=0, \\
\Phi _{0}^{\left( q\right) }=1, \\
\Phi _{n}^{\left( q\right) }=\Phi _{n-1}^{\left( q\right) }+\Phi
_{n-2}^{\left( q\right) }+\cdots +\Phi _{n-q-1}^{\left( q\right) }\text{ for
}n\geq 1.%
\end{array}%
\right.
\end{equation*}%
In \cite[2006]{bel1}, Belbachir and Bencherif proved that%
\begin{equation*}
\Phi _{n}^{\left( q-1\right) }=\sum_{k_{1}+2k_{2}+\cdots +qk_{q}=n}\binom{%
k_{1}+k_{2}+\cdots +k_{q}}{k_{1},k_{2},\cdots ,k_{q}},
\end{equation*}%
and, for $n\geq 1$%
\begin{equation*}
\Phi _{n}^{\left( q-1\right) }=\sum_{k=0}^{\left\lfloor n/\left( q+1\right)
\right\rfloor }\left( -1\right) ^{k}\frac{n-k\left( q-1\right) }{n-kq}\binom{%
n-kq}{k}2^{n-1-k\left( q+1\right) },
\end{equation*}%
leading to%
\begin{equation*}
\sum_{k_{1}+\cdots +qk_{q}=n}\binom{k_{1}+\cdots +k_{q}}{k_{1},\cdots ,k_{q}}%
=\sum_{k=0}^{\left\lfloor n/\left( q+1\right) \right\rfloor }\left(
-1\right) ^{k}\frac{n-k\left( q-1\right) }{n-kq}\binom{n-kq}{k}%
2^{n-1-k\left( q+1\right) }.
\end{equation*}

\noindent This is an analogous situation in writing above a multiple
summation with one symbol of summation. On the other hand, we establish a
connection between the ordinary multinomials and the generalized Fibonacci
sequence:

\begin{theorem}
We have the following identity%
\begin{equation}
\Phi _{n}^{\left( q\right) }=\sum_{l=0}^{qm-r}\binom{n-l}{l}_{q},
\label{uni}
\end{equation}%
where $m$ is given by the extended euclidean algorithm for division: $%
n=m\left( q+1\right) -r,$ $0\leq r\leq q.$
\end{theorem}

\begin{proof}
We have
\begin{eqnarray*}
\Phi _{n}^{\left( q\right) } &=&\sum_{k_{1}+2k_{2}+\cdots +\left( q+1\right)
k_{q+1}=n}\binom{k_{1}+k_{2}+\cdots +k_{q+1}}{k_{1},k_{2},\cdots ,k_{q+1}} \\
&=&\sum_{L\geq 0}\ \sum_{k_{1}+2k_{2}+\cdots +\left( q+1\right) k_{q+1}=n}%
\binom{L}{k_{1},k_{2},\cdots ,k_{q+1}} \\
&=&\sum_{L\geq 0}\ \sum_{k_{2}+2k_{3}+\cdots +qk_{q+1}=n-L}\binom{L}{%
k_{1},k_{2},\cdots ,k_{q+1}} \\
&=&\sum_{L\geq \frac{n}{q+1}}^{n}\binom{L}{n-L}_{q}.
\end{eqnarray*}

Now consider the unique writing of $n$ given by the extended euclidean
algorithm for division: $n=m\left( q+1\right) -r,$ $0\leq r<q+1$ $%
\rightarrow \frac{n}{q+1}=m-\frac{r}{q+1},$ which gives%
\begin{equation*}
\Phi _{n}^{\left( q\right) }=\sum_{k=0}^{qm-r}\binom{m+k}{qm-r-k}%
_{q}=\sum_{l=0}^{qm-r}\binom{m+k}{\left( q+1\right) k+r}_{q}=%
\sum_{l=0}^{qm-r}\binom{n-l}{l}_{q}.
\end{equation*}
\end{proof}

\noindent As an immediate consequence of Theorem 3, we obtain the following
identities%
\begin{eqnarray*}
\Phi _{\left( q+1\right) m}^{\left( q\right) } &=&\sum_{l=0}^{qm}\binom{%
\left( q+1\right) m-l}{l}_{q}=\sum_{k=0}^{qm}\binom{m+k}{\left( q+1\right) k}%
_{q}, \\
\Phi _{\left( q+1\right) m-1}^{\left( q\right) } &=&\sum_{l=0}^{qm-1}\binom{%
\left( q+1\right) m-l-1}{l}_{q}=\sum_{k=0}^{qm}\binom{m+k}{\left( q+1\right)
k+1}_{q}, \\
&&............ \\
\Phi _{\left( q+1\right) m-r}^{\left( q\right) } &=&\sum_{l=0}^{qm-r}\binom{%
\left( q+1\right) m-l-r}{l}_{q}=\sum_{k=0}^{qm}\binom{m+k}{\left( q+1\right)
k+r}_{q}.
\end{eqnarray*}

\noindent For $q=1,$ we find the classical Fibonacci sequence:%
\begin{equation*}
F_{-1}=0,\ F_{0}=1,\ F_{n+1}=F_{n}+F_{n-1},\ \text{for }n\geq 0.
\end{equation*}

\noindent Thus, we obtain the well known identity%
\begin{equation*}
F_{n}=\sum_{l=0}^{\left\lfloor \frac{n}{2}\right\rfloor }\binom{n-l}{l}.
\end{equation*}

\section{Exponential partial Bell partition polynomials}

\noindent In this section, we establish a connection of the ordinary
multinomials with exponential partial Bell partition polynomials $%
B_{n,L}\left( t_{1},t_{2},\ldots \right) $ which are defined (see Comtet
\cite[1970, p. 144]{com}) as follows%
\begin{equation}
\frac{1}{L!}\left( \sum_{m\geq 1}\frac{t_{m}}{m!}x^{m}\right)
^{L}=\sum_{n\geq L}B_{n,L}\frac{x^{n}}{n!},\ L=0,1,2,\ldots .  \label{bell}
\end{equation}

\noindent An exact expression of such polynomials is given by
\begin{equation*}
B_{n,L}\left( t_{1},t_{2},\ldots \right) =\sum_{\substack{ %
k_{1}+2k_{2}+\cdots =n  \\ k_{1}+k_{2}+\cdots =L}}\frac{n!}{%
k_{1}!k_{2}!\cdots \left( 1!\right) ^{k_{1}}\left( 2!\right) ^{k_{2}}\cdots }%
t_{1}^{k_{1}}t_{2}^{k_{2}}\cdots .
\end{equation*}

\noindent In this expression, the number of variables is finite according to
$k_{1}+2k_{2}+\cdots =n.$

\noindent Next, we give some particular values of $B_{n,L}:$%
\begin{eqnarray}
B_{n,L}\left( 1,1,1,\ldots \right) &=&\QDATOPD\{ \} {n}{L}\text{
Stirling
numbers of second kind,}  \notag \\
B_{n,L}\left( 0!,1!,2!,\ldots \right) &=&\QDATOPD[ ] {n}{L}\text{
Stirling
numbers of first kind,}  \notag \\
B_{n,L}\left( 1!,2!,3!,\ldots \right) &=&\frac{n!}{L!}\binom{n-1}{n-L}.
\label{ana}
\end{eqnarray}

\noindent In \cite[2005]{abb}, Abbas and Bouroubi give several extended
values of $B_{n,L}.$

\ \ \ \ \ \ \ \ \ \ \ \ \ \ \ \ \ \ \ \ \ \ \ \ \

\noindent The connection with ordinary multinomials is given by the
following result:

\begin{theorem}
We have the following identity%
\begin{equation}
B_{n,L}\left( 1!,2!,\ldots ,\left( q+1\right) !,0,\ldots \right) =\frac{n!}{%
L!}\binom{L}{n-L}_{q}.  \label{anaa}
\end{equation}
\end{theorem}

\begin{proof}
Taking in (\ref{bell}) $t_{m}=m!$ for $1\leq m\leq q+1$ and zero otherwise,
we obtain%
\begin{equation*}
\left( x+\cdots +x^{q+1}\right) ^{L}=L!\sum_{n-L\geq 0}B_{n,L}\left(
1!,2!,\ldots ,\left( q+1\right) !,0,\ldots \right) \frac{x^{n}}{n!},
\end{equation*}%
from which it follows%
\begin{equation*}
\sum_{a\geq 0}\binom{L}{a}_{q}x^{a}=\sum_{n-L\geq 0}\frac{L!}{n!}%
B_{n,L}\left( 1!,2!,\ldots ,\left( q+1\right) !,0,\ldots \right) x^{n-L}.
\end{equation*}
\end{proof}

\begin{corollary}
Let $q\geq 1$, $L\geq 0$ be integers, and $a\in \left\{ 0,1,\ldots
,qL\right\} .$ For $q\geq a,$ we have the following identity%
\begin{equation*}
\binom{L}{a}_{q}=\binom{L+a-1}{a}.
\end{equation*}
\end{corollary}

\begin{proof}
Using the fact that $B_{n,L}\left( 1!,2!,\ldots ,\left( q+1\right)
!,0,\ldots \right) =B_{n,L}\left( 1!,2!,3!,\ldots \right) $ for $q+1\geq
n-L+1,$ we obtain $\binom{L}{n-L}_{q}=\binom{n-1}{n-L}$ for $q\geq n-L.$ We
conclude with $a=n-L.$
\end{proof}

\section{Convolution powers of discrete uniform distribution}

\noindent This section gives a connection between the ordinary multinomials
and the convolution power of the discrete uniform distribution. The right
hand side of identity (\ref{6b}) is a very well known expression. Indeed for
$q,L\in
\mathbb{N}
,$ let us denote by $U_{q}$ the $L^{th}$ convolution power of the discrete
uniform distribution%
\begin{equation*}
U_{q}:=\frac{1}{q+1}\left( \delta _{0}+\delta _{1}+\cdots +\delta
_{q}\right) \ \ \ \ \text{(}\delta _{a}\text{ is the Dirac measure),}
\end{equation*}%
then for $a\in
\mathbb{N}
$ (see de Moivre \cite[1731]{moi1} or \cite[1998]{hal}), with respect to the
counting measure, its density is given by%
\begin{equation}
P\left( U_{q}^{\star L}=a\right) =\frac{1}{\left( q+1\right) ^{L}}%
\sum_{j=0}^{\left\lfloor a/\left( q+1\right) \right\rfloor }\left( -1\right)
^{j}\binom{L}{j}\binom{a+L-\left( q+1\right) j-1}{L-1}.  \label{sss}
\end{equation}

\noindent Combining Theorem 1 and relation (\ref{sss}), we have the
following result:

\begin{corollary}
Using the above notations, we obtain the following identity%
\begin{equation*}
P\left( U_{q}^{\star L}=a\right) =\frac{\binom{L}{a}_{q}}{\left( q+1\right)
^{L}}.
\end{equation*}
\end{corollary}

\noindent It should be noted that the multinomials may be seen as the number
of favorable cases to the realization of the elementary event $\left\{
a\right\} .$

\noindent It is easy to show that the distribution of $U_{q}^{\star L}$ is
symmetric by relation (\ref{3}).

\begin{corollary}
We have the following identity%
\begin{eqnarray*}
\sum_{k=0}^{qL}k\binom{L}{k}_{q} &=&\left( q+1\right) ^{L}\frac{qL}{2}, \\
\sum_{k=0}^{qL}k_{q}^{2}\binom{L}{k} &=&\left( q+1\right) ^{L}\frac{qL}{2}%
\left( \frac{qL}{2}+\frac{q+2}{6}\right) , \\
\sum_{k=0}^{qL}k_{q}^{3}\binom{L}{k} &=&\left( q+1\right) ^{L}\left( \frac{qL%
}{2}\right) ^{2}\left( \frac{qL}{2}+\frac{q+2}{2}\right) ,
\end{eqnarray*}

More generally, for $m\geq 1,$ the following identity holds%
\begin{equation*}
\sum_{k=0}^{qL}k_{q}^{m}\binom{L}{k}=\left( q+1\right)
^{L}\sum_{i_{1}+i_{2}+\cdots +i_{L}=m}\binom{m}{i_{1},i_{2},\cdots ,i_{L}}%
u_{i_{1}}u_{i_{2}}\cdots u_{i_{L}},
\end{equation*}%
where $u_{i}$ is the $i$-th moment of the random variable $U_{q}.$
\end{corollary}

\begin{proof}
It suffices to compute the expectation of $U_{q}^{\star L}$ using, first the
density distribution and second the summation of uniform distributions.
\end{proof}

\begin{acknowledgement}
We are grateful to Professor Miloud Mihoubi for pointing our attention to
Bell polynomials.
\end{acknowledgement}

\end{document}